%% file: artikel4.tex
\documentclass[a4paper,11pt,english]{article}
\usepackage{amssymb}
\usepackage{mathtools}    
\usepackage[all]{xy}
\usepackage{amsthm}
\usepackage{mathdots}
\usepackage{amsmath}
\usepackage{verbatim}
\usepackage{pifont}
\usepackage{MnSymbol}
\usepackage{mathrsfs}
\usepackage[ansinew]{inputenc}
\usepackage{microtype}
\usepackage{hyperref}
\newtheorem{thm}{Theorem}[section]
\newtheorem{lemma}[thm]{Lemma}
\newtheorem{defi}[thm]{Definition}

\newtheorem{cor}[thm]{Corollary}

\newtheorem{prop}[thm]{Proposition}

\input{Befehle.tex}

\begin{document}

\thispagestyle{empty}
\begin{center}
\huge
\textsf{K-theory of hermitian symmetric spaces and root lattices}\\
\vspace*{1cm} \large
Dennis Bohle, Wend Werner\\
Fachbereich Mathematik und Informatik\\
Westf\"alische Wilhelms-Universit\"at\\
Einsteinstra\ss e 62\\
48149 M\"unster \\
\vspace*{0.5cm}
 dennis.bohle@math.uni-muenster.de,\\ wwerner@math.uni-muenster.de\\
\vspace*{0.3cm}
\end{center}

\begin{abstract}
This is a companion to a recent investigation of K-theoretical invariants for symmetric spaces.
We introduce a new class of cycles in K-groups, which are connected to elements of an underlying root lattice.
This will be needed for a K-theoretical classification of inductive limits.
\end{abstract}

\section{Introduction}

In a parallel investigation \cite{BoWe4}, inductive limits of bounded symmetric domains are classified
through a variant of K-theory. The latter has recently been shown to yield a homological method for the
classification of the classical bounded symmetric domains \cite{BoWe2}.

In finite dimensions, the Cartan/Serre approach requires a Hilbert space structure on the vector spaces containing
the root system of the underlying Lie algebra. Also in \cite{BoWe2} roots (or, more precisely, grids) have been used to
mark K-groups in order to achieve a complete invariant. Using root systems is problematic if treating inductive limits:
The Hilbertian structure, along the limit, disappears completely, and the connecting
morphisms may be of higher multiplicity, which makes it difficult to keep track of the underlying root systems.

In the sequel we try to promote the idea to use a different invariant, a marked subset of the group K$_0$,
which is closely linked to root lattices.
As can be seen in \cite{BoWe4}, this new ingredient is very helpful in infinite dimensions.

We will provide some background information in the next section, define and motivate the new K-theoretical invariant in the third, and then set out for a calculation in the final section.

\section{Partial isometries, Grids and Roots}
Important for the following is the Kaup-Koecher-Loos approach to symmetric complex domains. It features
a threefold product, defined on the Banach space $E$
surrounding the domain $U$. This product turns $E$ into an algebraic object, a JB*-triple system, more general
than a C*-algebra.
\begin{defi}
A Banach space $Z$ together with a trilinear, continuous and w.r.t. the outer
variables symmetric mapping $\{\cdot,\cdot,\cdot\}:Z^3\to Z$ is called a \emph{JB*-triple}, iff
\begin{itemize}
\item[(a)] $||\{x,x,x\}||=||x||^3$ for all $x\in Z$,
\item[(b)] With the operator $x\Box y$ defined on $Z$ by $(x\Box y)(z)=\{x,y,z\}$,
$ix\Box y$ is a derivation,
\item[(c)] $x\Box x$ has non-negative spectrum, and $\exp(it(x\Box x))$ is a 1-parameter group of isometries.
\end{itemize}
\end{defi}
It turns out that complex bounded domains in Banach spaces are, up to biholomorphic bijections, just the open unit
balls of JB*-triple systems \cite{Kaup-AlgebraiccharacterizationofsymmetriccomplexBanachmanifolds}. In finite dimension,
a Wedderburn type result shows that JB*-triple systems can be decomposed into direct sums of \emph{Cartan factors},
which will be detailed below.

The structure of a JB*-triple system is intimately connected to the Lie algebra of all complete holomorphic vector fields
defined on its open unit ball. (Recall that a \emph{complete vector field} is supposed to have a flow defined for all times.)
A variant of this relationship is the Tits-Kantor-Koecher Algebra of $Z$, which is given as the Lie (sub)algebra of quadratic vector fields
\begin{equation*}
    \ftkk(Z)=\set{z\mapsto a+\sum_i\{v_i,w_i,z\}+\{z,b,z\}}{a,b,v_i,w_i\in Z}.
\end{equation*}
It is 3-graded, $\ftkk(Z)=\ftkk(Z)_{-1}\oplus\ftkk(Z)_0\oplus\ftkk(Z)_1$,
where $\ftkk(Z)_0$ is the Lie subalgebra of linear vector fields, and
\begin{equation*}
    \ftkk(Z)_\pm=\set{z\mapsto a\pm\{z,a,z\}}{a\in Z}.
\end{equation*}

Root systems of $\ftkk(Z)$ can be recovered from the JB*-triple structure as they leave an imprint on the set
of \emph{tripotents}, i.e.\ elements $z\in Z$ with $\{z,z,z\}=z$. Central to this relationship is the concept of a
\emph{grid} or, more generally of a \emph{cog}, both sets of tripotens satisfying a number of relations.
This connection has been thoroughly explored by Neher \cite{Neher-3gradedrootsystemsandgridsinJordantriplesystems,
Neher-Systemesderacines3gradues} who has shown how to use grids for a classification of finite dimensional bounded symmetric
domains based on generators and relations.

We recall the definition of a root system.
Let $X$ be a real, finite-dimensional vector space with scalar product $\langle\cdot,\cdot\rangle$. A subset $R\sub X$ is called a \emph{root system} iff
\begin{description}
\item[(a)] $R$ is finite, generates $X$ and does not contain $0$.
\item[(b)] For every root $\alpha\in R$ we have $s_\alpha(R)=R$, where $s_\alpha(x)=x-2\frac{\langle x,\alpha\rangle}{\langle\alpha,\alpha\rangle}\alpha$ is the reflection in $\alpha$.
\item[(c)] $(\alpha,\beta):=2\frac{\langle\alpha,\beta\rangle}{\langle\beta,\beta\rangle}\in\mathbb{Z}$ for all $\alpha,\beta \in R$.
\item[(d)] For every $\alpha\in R$, we have $\R\alpha\cap R=\{\pm\alpha\}$.
\end{description}
A root system $R\neq \emptyset$ is called \emph{irreducible} if $R$ cannot be decomposed into two orthogonal non-empty subsets.
The root system $R$ is \emph{3-graded} if additionally there exist $R_{-1},R_0,R_1\sub R$ with
\begin{description}
\item[(d)] $R=R_1\dot\cup R_0\dot\cup R_{-1}$ (disjoint union).
\item[(e)] $R_{-1}=-R_1$.
\item[(f)] $R_0=\{\alpha-\beta;\alpha,\beta\in R_1, \alpha \neq \beta, \langle\alpha,\beta\rangle=0\}$.
\item[(g)] If $\alpha,\beta\in R_1$ then $\alpha+\beta\not\in R$.
\item[(h)] If $\alpha\in R_0,\beta\in R_1$ and $\alpha+\beta\in R$ then $\alpha+\beta\in R_1$.
\end{description}
Due to the conditions (e) and (f) the grading of a root system is completely determined by its $(R_1)$ part.

It can be shown that for each (closed) cog in $Z$ there is a 3-graded root system
$R$ for $\ftkk(Z)$ as well as a bijection respecting the relations among cog and,
respectively, root elements. Irreducible 3-graded root systems correspond to \emph{connected grids}.

The relation between roots and grids is as follows:
Identify $Z$ with the 1-part of its Kantor-Koecher-Tits algebra $\ftkk(Z)$. The root system $R$ of
$\ftkk_0(Z)$ with respect to a Cartan subalgebra $\fh\sub\ftkk_0( V)$ is then naturally
3-graded, and the spaces $\ftkk_i(Z)$, $i=\pm 1,0$, are spanned by the root
spaces $\ftkk(Z)_\alpha$ with $\alpha\in R_i$. Moreover, relative to a properly chosen
$\fh\sub\ftkk_0( V)$, the root spaces $\ftkk(Z)_\alpha$, $\alpha\in R_i$, are exactly the one dimensional
subspaces spanned by the elements of the grid.

Every connected grid $\cG$ is associated to one of the following standard grids.
We include a list of the $3$-graded root systems corresponding to the standard grids
and thus to the classical Cartan factors. For details see \cite{Neher-Jordantriplesystemsbythegridapproach}, \cite[§3]{NeherErhard-Liealgebrasgradedby3gradedrootsystemsandJordanpairscoveredbygrids}.

\paragraph{Cartan factor $\mathbf{C^1_{n,m}}$, of type I}
This is the space of complex $n\times m$-matrices $\M_{n,m}(\mathbb{C})$, the \emph{rectangular factor},
with triple product $1/2(AB^*C+CB^*A)$.
Its (reduced) grid $\mathcal{R}(n,m)$, $n,m\geq 1$ is called rectangular and is given by the usual matrix units
$\set{E_{i,j}}{1\leq i\leq m,\ 1\leq j\leq n}$. The corresponding root system $R$ lies in the subspace
\begin{equation*}
X=\set{\sum_{1\leq k\leq m+n} s_ke_k\in\ell^2(m+n)}{\sum_{1\leq k\leq m+n} s_k=0}
\end{equation*}
of  the (real)
Hilbert space $\ell^2(m+n)$. (Here, as in the following, $e_i$ are elements of an orthonormal basis.)
It is given by
\begin{equation*}
R=\mathscr{A}_{m+n-1}=\set{e_k-e_l}{1\leq k,l\leq m+n,k\neq l}
\end{equation*}
with grading induced by
\begin{equation*}
R_1=\set{e_i-e_j}{1\leq i\leq m,\ 1\leq j\leq n}.
\end{equation*}

\paragraph{Cartan factor $\mathbf{C^2_n}$, of type II}
A space of skew-symmetric, complex $n\times n$-matrices, $n\geq 4$, called \emph{symplectic factor}.
The standard example of the symplectic grid $\cS(n)$, $n\geq 5$, is $\set{E_{i,j}-E_{j,i}}{1\leq i,j\leq n, i<j}$,
the corresponding $3$-graded root system is a subset of $\ell^2(n)$, and is given by
$$
R=\mathscr{D}_n=\set{\pm e_i\pm e_j}{1\leq i,j\leq n, i\neq j},
$$
with grading given by
$$
R_1=\set{e_i+e_j}{1\leq i,j\leq n, i\neq j}.
$$

\paragraph{Cartan factor $\mathbf{C^3_n}$, of type III}
This is the \emph{hermitian factor}, consisting of symmetric complex $n\times n$-matrices, $n\geq 2$.
For the hermitian grid $\cH(n)$, $n\geq 2$, the standard example is
$\set{E_{i,j}+E_{j,i}}{1\leq i,j\leq n, i\neq j}\cup\set{E_{i,i}}{1\leq i\leq n}$ the associated root system is
$$
R=\mathscr{C}_n=\set{\pm 2e_i}{1\leq i\leq n}\cup\set{\pm e_i\pm e_j}{1\leq i,j\leq n}
$$
with $1$-part
$$
R_1=\set{e_i+e_j}{1\leq i,j\leq n},
$$
again as a subset of $\ell^2(n)$.

\paragraph{Cartan factor $\mathbf{C^4_n}$, of type IV}
The $n+1$-dimensional \emph{spin factor}, $n\geq 2$ is the
closed linear span of selfadjoint matrices $1,s_1,\ldots,s_n$, satisfying
\begin{equation*}
s_is_j+s_js_i=2\delta_{ij}1
\end{equation*}
for all $i,j\in\{1,\ldots,n\}$. Such a set of matrices is called
a \emph{spin system}\index{spin system}. The standard way of constructing one
is to start with the Pauli matrices
\begin{equation*}
\sigma_1:=\begin{pmatrix}
                    1 & 0 \\
                    0 & -1
            \end{pmatrix}\quad
\sigma_2:=\begin{pmatrix}
                    0 & 1 \\
                    1 & 0
            \end{pmatrix}\quad
\sigma_3:=\begin{pmatrix}
                    0 & i \\
                   -i & 0
            \end{pmatrix}.
\end{equation*}
For a matrix $\sigma$, write $\sigma^k$ to denote its k-fold tensor product. Then
\begin{align*}
s_0&=id^{n},\\
s_1&=\sigma_1\otimes\id^{n-1},\\
s_2&=\sigma_2\otimes\id^{n-1},\\
s_3&=\sigma_3\otimes\sigma_1\otimes\id^{n-2},\\
s_4&=\sigma_3\otimes\sigma_2\otimes\id^{n-2},\\
s_{2l+1}&=\sigma_3^{l}\otimes\sigma_1\otimes\id^{n-l-1},\\
s_{2l+2}&=\sigma_3^{l}\otimes\sigma_2\otimes\id^{n-l-1},
\end{align*}
$1\leq l\leq n-1$, becomes a spin system which generates the odd dimensional spin factor.
If we drop the last idempotent $s_{2n}$ we get a spin system which generates an even dimensional spin factor.
The spin grid $\Sp(n)=\set{u_i,\widetilde{u}_i}{i=1,\ldots,n}$, $n\geq 1$ is obtained in the following way.
(We follow here an idea from \cite{Bohle}, which in turn is based on
\cite{nealrusso-Contractiveprojectionsandoperatorspaces}.)
If $n$ is odd (i.e.\
the corresponding spin factor is of even dimension) the elements of the $\Sp(n)$
are given by
\begin{gather*}
u_1=\frac{1}{2}\lr{\id-s_1}\quad \widetilde{u}_1=-\frac{1}{2}\lr{\id+s_1}\quad\text{and}\\
u_{k+1}=\frac{1}{2}\lr{s_{2k}+is_{2k+1}},\quad\widetilde{u}_{k+1}=\frac{1}{2}\lr{s_{2k}-is_{2k+1}}\quad
\text{for $k=1,\ldots,\frac{1}{2}\lr{n-1}$},
\end{gather*}
whereas in the case that $n$ is even,
$\Sp(n)=\set{u_i,\widetilde{u}_i}{i=1,\ldots,n}\cup\{u_0\}$ with
$u_0=s_n$.
The $3$-graded root system associated to $\Sp(n)$, $n$ even, is a root system of type
$\mathscr{D}_{n+1}$ in the vector space $\ell^2(n)\oplus\mathbb{R}e_\infty$,
$$
R=\mathscr{D}_{n+1}=\set{\pm e_i\pm e_j}{i,j=1,\ldots,n,\infty, i\neq j},
$$
where the grading comes from
$$
R_1=\set{e_{\infty}\pm e_i}{i=1,\ldots,n}.
$$
For a spin grid $\Sp(n)$ with $|\Sp(n)|$ odd the associated root system $R$ is of type $\mathscr{B}_{\text{\tiny $n+1$}}$,
$$
R=\mathscr{B}_{\text{\tiny $n+1$}}=\set{\pm e_i\pm e_j}{i,j=1,\ldots,n,\infty, i\neq j}
\cup\set{\pm e_i}{i=1,\ldots,n,\infty},
$$
$3$-graded with 1-part
\[
R_1=\set{e_\infty\pm e_i}{i=1,\ldots,n}\cup\{e_\infty\}.
\]

\paragraph{Exceptional factors}
There are two exceptional factors in dimensions $16$ and $27$.
And the two exceptional grids
\begin{description}
\item[(V)] Bi-Caley grid.
\item[(VI)] Albert grid.
\end{description}
Unfortunately, they do not play a role here.

These factors are mutually non-isomorphic with the exception
of $\M_{2,2}$, the symplectic factor for $n=4$, and the hermitian factor for $n=2$, which are spin factors,
of dimensions 4, 6 and 3, respectively.

\section{The K-JB* invariant}
Recall the following definitions from \cite{BoWe2}. For unexplained notation on ternary rings of operators
(TROs, for short), we refer the reader to \cite{BlecherLeMerdy-Operatoralgebrasandtheirmodules}.

K-groups for JB*-triples are obtained in two steps. The first defines
K-groups for TROs: For a TRO $T$, $K_*(T)$ is the K-group of its left C*-algebra.
As each TRO-morphism $\psi$ yields
C*-morphism $\mathcal{L}(\psi)$ and $\mathcal{R}(\psi)$ between both, left and right C*-algebras,
we may define $K_*(\psi)$ in a functorial way.
Note that for TROs, we stick to writing $K_*(\cdot)$ for their K-groups, as these groups coincide with their C*-counterparts
should they themselves be C*-algebras.
The second step involves the enveloping
TRO of a JB*-triple \cite{BoWe1,BunceFeelyTimoneyI}. To each JB*-triple corresponds a
\emph{universal (enveloping) TRO} $T^*(Z)$, as well as a canonical embedding $\rho_Z:Z\to T^*(Z)$ such
that the image $\rho_Z(Z)$ generates $T^*(Z)$ as a TRO.
Furthermore, each JB*-morphism uniquely lifts to a TRO-morphism between the respective universal TROs,
and the emerging functor, denoted $\tau$, has all properties needed for the ensuing K-theory.
For a JB*-triple $Z$ and a JB*-morphism $\phi:Z\to W$, we define
\begin{equation*}
    K_*^{\mathrm{JB}^*}(Z)=K_*^{\mathrm{TRO}}(\tau(Z)),
\end{equation*}
as well as
\begin{equation*}
    K_*^{\mathrm{JB}^*}(\phi)=K_*^{\mathrm{TRO}}(\tau(\phi))
\end{equation*}
This K-functor has the usual properties that one would expect from it,
except stability, which already has a bad start, as it is not possible,
in general, to equip the space of matrices with entries from a JB*-triple
with the structure of a JB*-triple in a natural way.

A complete isomorphism invariant for finite dimensional JB*-triples is obtained by
equipping K-groups with further structure, consisting of some distinguished subsets of $K_0(Z)$. The first is
the semigroup $K_0(Z)^+$, consisting of all classes
of projections themselves, that is, the set obtained before, in the final step, the Grothendieck construction
is applied in order to produce $K_0(Z)$. The next is the scale of the left C*-algebra of $T^*(Z)$, i.e.\
classes of projections in $\mathcal{L}(T^*(Z))$. We access the scale of the right C*-algebra of $T^*(Z)$ with the
aid of a canonical mapping which (at least in case of separable C*-algebras) is defined in the following way.
It can be shown that the canonical mappings $\iota_\cL:\cL(T)\to\Li(T)$ and $\iota_\cR:\cR(T)\to\Li(T)$ induce isomorphisms
\begin{equation*}
K_0(\iota_\cR):K_0(T)\to K_0(\Li(T))\quad\text{and}\quad
K_0(\iota_\cR):K_0(\cR(T))\to K_0(\Li(T)).
\end{equation*}
We then let the image of $K_0(\iota_\cR)K_0(\iota_\cR)^{-1}$ represent the scale of $\mathcal{R}(T^*(Z))$
in $K_0(\cR(T))$.

So far, the data we have assembled wouldn't suffice for a complete isomorphism invariant. In \cite{BoWe2}, K-classes generated
by grid elements were shown to provide just the right amount of additional information missing before. These classes,
however, show an inappropriate behavior when dealing with morphisms other than automorphisms. This defect may be
overcome by considering the lattice generated by grid elements. If, however, we follow along the lines of the correspondence
between roots and grid elements, a more practicable approach comes into view. Consider an element
\begin{equation*}
    \beta=\sum k_i\alpha_i,\qquad k_i\in\Z.
\end{equation*}
of the root lattice
of $\ftkk(Z)$, generated by a root system for a properly chosen Cartan-subalgebra $\fh\sub\ftkk(Z)_0$. Then all
elements $\sum\lambda_i g_i$, $\lambda_i\in\C$ and $g_i\in\ftkk(Z)_{\alpha_i}$ belong to
\begin{equation*}
    \ftkk(Z)_\beta=\set{z\in Z}{\text{$[h,z]=\beta(h)z$ for all $h\in\fh$}}.
\end{equation*}
But, as all tripotents of $Z$ can be written in such a way, we use
\begin{equation*}
    \Delta(Z):=\set{\left[\rho_Z(g)\rho_Z(g)^*\right]\in K_0^{\text{\tiny JB*}}(Z)}
    {g\in \Tri(Z)}
\end{equation*}
as an additional ingredient for the isomorphism invariant. This set carries as much information as the root
lattice of $\ftkk(Z)$ might provide on the level of K-groups. Summing up,
\begin{defi}
Let $Z$ be a $JB^*$-triple system. The
K-JB* invariant\index{K-grid invariant of a $JB^*$-triple system} of $Z$ is the tuple
\[
\mathcal{KJB^*}(Z):=\left(K_0^{\text{\tiny JB*}}(Z),K_0^{\text{\tiny JB*}}(Z)_+,\Sigma^{\text{\tiny JB*}}
_\mathcal{L}(Z),\Sigma^{\text{\tiny JB*}}_\mathcal{R}(Z),\Delta(Z)\right),
\]
where $K_0^{\text{\tiny JB*}}(Z)_+:=K_0(T^*(Z))_+$, $\Sigma_\cL^{\text{\tiny JB*}}(Z)$ and
$\Sigma^{\text{\tiny JB*}}_\cR(Z)$ are the left and right scale of the TRO $T^*(Z)$ and
$\Delta(Z)$ is the set of equivalence classes coming from the set of tripotents in $Z$.
\end{defi}

\section{Calculation}
Since in finite dimensional C*-algebras projections are equivalent iff they are unitarily equivalent, the following
lemma follows from \cite[Lemma 3.1]{BoWe4}.
\begin{lemma}\label{tripot-equiv}
Let $u,v$ be tripotents in a finite dimensional JB*-triple $Z$. Then the following are equivalent.
\begin{description}
\item[(i)]
$u$ and $v$ yield the same class in $K_0^{\text{\tiny JB*}}(Z)$
\item[(ii)]
There is a unitary $U\in\cL(T^*(Z))$ such that $Uv=u$.
\item[(iii)]
There is an automorphism of $Z$ mapping $u$ onto $v$.
\end{description}
\end{lemma}
The TRO $T^*(Z_1)$ is the finite sum of rectangular matrix algebras $T^*(Z_1)\simeq\bigoplus_{i=1}^p\M_{n_i,m_i}$
\cite{Smith-Finitedimensionalinjectiveoperatorspaces}. This structure is perfectly reflected by the double-scaled
ordered $K_0$-group of $T^*(Z_1)$.

\begin{defi}
Let $\varphi:K_0^{\text{\tiny JB*}}(Z_1)\to K_0^{\text{\tiny JB*}}(Z_2)$ be a map. We say that $\varphi$ is a K-JB* morphism if $\varphi$ is a group morphism with
\begin{gather*}
\vp(K_0^{\text{\tiny JB*}}(Z)_+)\sub K_0^{\text{\tiny JB*}}(W)_+,\qquad
\varphi(\Sigma^{\text{\tiny JB*}}_\mathcal{L}(Z_1))\sub\Sigma^{\text{\tiny JB*}}_\mathcal{L}(Z_2),\\
\varphi(\Sigma^{\text{\tiny JB*}}_\mathcal{R}(Z_1))\sub\Sigma^{\text{\tiny JB*}}_\mathcal{R}(Z_2),
\quad\text{and}\quad
\varphi(\Delta(Z_1))\sub\Delta(Z_2).
\end{gather*}
\end{defi}

The proof of the following result is obvious (and follows along the lines of \cite[Proposition 4.5]{BoWe2})

\begin{prop}\label{kgrid inv additiv}
Let $Z_1$ and $Z_2$ be finite dimensional JB*-triple systems, then there exists a $\mathcal{KJB^*}$-isomorphism of $K_0$-groups
\begin{equation*}
\mathcal{KJB^*}(Z_1\oplus Z_2)\simeq \mathcal{KJB^*}(Z_1)\oplus \mathcal{KJB^*}(Z_2).
\end{equation*}
\end{prop}

\begin{prop}
Let $Z$ be a $JB^*$-triple system
\begin{itemize}
\item[(a)] Suppose $Z$ is the finite-dimensional type I Cartan factor $Z=C^1_{n,m}$.
\begin{itemize}
\item[(i)]
If $Z$ is isometric to a finite-dimensional Hilbert space, i.e.\ $Z=C^1_{1,n}$, $n\in \N$, then $\mathcal{KJB^*}(C^1_{1,n})$ is given by
\begin{align*}
K_0^{\text{\tiny JB*}}(C^1_{1,n})\;\;&=\Z^n,\\
K_0^{\text{\tiny JB*}}(C^1_{1,n})_+&=\N_0^n,\\
\Sigma_\mathcal{L}^{\text{\tiny JB*}}(C^1_{1,n})\;\;&=\prod_{k=1}^n\left\{1,\ldots,\begin{pmatrix}n\\
k\end{pmatrix}\right\},\\
\Sigma_\mathcal{R}^{\text{\tiny JB*}}(C^1_{1,n})\;\;&=\prod_{k=1}^n\left\{1,\ldots,\begin{pmatrix}n\\ k-1\end{pmatrix}\right\}\text{ and}\\
\Delta(C^1_{1,n})\;\;&=\left\{\left(\begin{pmatrix} n-1 \\ 0\end{pmatrix},\begin{pmatrix} n-1 \\ 1\end{pmatrix},\ldots,\begin{pmatrix} n-1 \\ n-1\end{pmatrix}\right)\right\}.
\end{align*}

\item[(ii)]  If $n,m\geq 2$, then $\mathcal{KJB^*}(C^1_{n,m})$ is given by
\begin{align*}
\mathcal{KJB^*}(C^1_{n,m})=(\Z^2,\N_0^2,&\{1,\ldots,n\}\times\{1,\ldots,m\},\\
&\{1,\ldots,m\}\times\{1,\ldots,n\},\\
&\{1,\ldots,n\}\times\{1,\ldots,m\}).
\end{align*}
\end{itemize}
\item[(b)] Let $Z$ be isometric to a Cartan factor of type II with $\dim Z\geq 10$. Then  \[
\mathcal{KJB^*}(C^2_{n})=(\Z,\N_0,\{1,\ldots,n\},\{1,\ldots,n\},\{2,4\ldots,k\}),
\]
where $k$ is the greatest even integer less or equal to $n$.
\item[(c)] If  $Z$ is $JB^*$-triple isomorphic to the finite-dimensional Cartan factor $C^3_n$, then \[
\mathcal{KJB^*}(C^3_n)=\left(\Z,\N_0,\{1,\ldots,n\},\{1,\ldots,n\},\{1,\ldots,n\}\right).
\]
\item[(d)]  Let $Z$ be a finite-dimensional spin factor with $\dim Z=k+1.$
\begin{itemize}

 \item[(i)] If $Z$ is of even dimension, i.e.\ $k=2n-1$, $n\geq 2$, then the K-$JB^*$ invariant of $Z$ is given by
\begin{align*}
K_0^{\text{\tiny JB*}}(Z)\;\;&=\Z^2,\\
K_0^{\text{\tiny JB*}}(Z)_+&=\N_0^2,\\
\Sigma^{\text{\tiny JB*}}_\cL\left(Z\right)\;\;&=\left\{1,\ldots,2^{n-1}\right\}^2,\\
\Sigma^{\text{\tiny JB*}}_\cR\left(Z\right)\;\;&=\left\{1,\ldots,2^{n-1}\right\}^2,\\
\Delta\left(Z\right)\;\;&=\left\{(2^{n-2},2^{n-2}),(2^{n-1},2^{n-1})\right\}.
\end{align*}

\item[(ii)] If $Z$ is of odd dimensions, i.e.\ $k=2n$, $n\geq 2$, then $\mathcal{KJB^*}(Z)$ has the components
\begin{align*}
K_0^{\text{\tiny JB*}}\left(Z\right)\;\;&=\Z,\\
K_0^{\text{\tiny JB*}}(Z)_+&=\N_0,\\
\Sigma_\cL^{\text{\tiny JB*}}\left(Z\right)\;\;&=\{1,\ldots,2^{n}\},\\
\Sigma_\cR^{\text{\tiny JB*}}\left(Z\right)\;\;&=\{1,\ldots,2^{n}\},\\
\Delta\left(Z\right)\;\;&=\left\{2^{n-1},2^{n}\right\}.
\end{align*}
\end{itemize}
\end{itemize}
\begin{proof}

All the $K_0$-groups, their positive cones and their scales were computed in \cite{BoWe2}, so we have to compute the new invariant $\Delta(Z)$ for all Cartan factors $Z$. In light of Lemma~\ref{tripot-equiv}, the result for the (non-Hilbertian)
rectangular, symplectic and hermitian cases are rather obvious.
We prove the remaining cases, Hilbert spaces and spin factors.

The tripotents $e$ in a Hilbert space, when
the latter is viewed as a Cartan factor of type I (and $\rank$ 1), are the norm-one elements. Thus, whenever $e,f$ are
tripotents, there exists an automorphism $U$, i.e.\ a unitary map, mapping $e$ to $f$. Hence, $\Delta(Z)$ coincides with the classes generated by the grid elements which in light of \cite[Proposition 5.2.]{BoWe2} completes the proof of (a)(i).

If $Z$ is a spin factor of dimension $n$, tripotents belong to one of two distinct classes,
which can be seen as follows.
Recall that on a spin factor $Z=\mathbb{C}^n$ there is always a scalar product as well as an involution
$a\mapsto\bar{a}$ so that
\begin{equation*}
\{a,b,c\}=\langle a,b\rangle c+\langle c,b\rangle a-\langle a,\bar{c}\rangle \bar{b}
\end{equation*}
for $a,b,c\in\mathbb{C}^n$. Thus $e\in Z$ is a tripotent iff
$$
e=\{e,e,e\}=2\langle e,e\rangle e-\langle e,\bar{e}\rangle \bar{e},
$$
and we have the following possibilities:
\begin{description}
\item[(a)] $\langle e,\bar{e}\rangle =0$ and $\langle e,e\rangle =\fr{2}$.
\item[(b)] $\langle e,\bar{e}\rangle \neq0$, $e=\mu r$ with $\mu\in\mathbb{C}$, $|\mu|=1$ and $r\in \mathbb{C}^n$, $\bar{r}=r$ and $\langle r,r\rangle =1$.
\end{description}
Following the notation of \cite[Chapter 3]{Friedman-Physicalapplicationsofhomogeneousballs}, the tripotents in (a) are called minimal and those in (b) maximal. We will show that that two minimal tripotents yield equivalent elements in $K_0(Z)$ and likewise for two maximal tripotents.

Let $e$ and $f$ be maximal tripotents. Thus $e=\mu_1r_1$ and $f=\mu_2r_2$ as above. It is well known \cite[Corollary 5]{Harris-Boundedsymmetrichomogeneousdomainsininfinitedimensionalspaces} that a linear map $\Phi$ is an automorphism of
$Z$ iff $\Phi=\lambda\Phi_0$, where $\lambda$ is a complex scalar of modulus $1$ and $\Phi_0$ is unitary (for the scalar
product $\langle \cdot,\cdot\rangle $) such that $\Phi_0(\bar{x})=\overline{\Phi_0(x)}$ for all $x\in Z$.
Let $\Phi_{00}$ be an orthogonal mapping defined on the self-adjoint part of $Z$, mapping $r_1$ to $r_2$.
Denote by $\Phi_0$ its complexification (which then is an automorphism of the spin factor).
If we define $\Phi:=\bar{\mu_1}\mu_2\Phi_0$, then $\Phi$ is an automorphism of
$Z$ mapping $e$ to $f$ and the corresponding projections are equivalent in $K_0(Z)$.

It follows from \cite[3.1.4]{Friedman-Physicalapplicationsofhomogeneousballs} that if we decompose a minimal tripoten $v$ as
\begin{equation*}
v=x+iy,\qquad\text{$x,y$ self-adjoint},
\end{equation*}
then
$$
|x|=|y|=\fr{2}\quad\text{and}\quad\langle x,y\rangle=0.
$$
Let now $e=x_1+iy_1$ and $f=x_2+iy_2$ be minimal tripotents decomposed in this way.
Pick $e_1^i,\ldots,e_{n-2}^i$, such that $2x_i,2y_i,e_1^i,\ldots,e_{n-2}^i$ is an orthonormal basis of
the self-adjoint part of $Z$, for $i=1,2$. Let $\Phi_{00}$ be the linear mapping with $\Phi_{00}(x_1)=x_2$,
$\Phi_{00}(y_1)=y_2$ and $\Phi_{00}(e^1_k)=e^2_k$ for $k=1,\ldots,n-2$. The complexification $\Phi_0$ of
$\Phi_{00}$ is obviously an automorphism of $Z$ that maps $e$ to $f$ and therefore the elements in $K_0(Z)$
that correspond to $e$ and $f$ are equivalent.

Along the same lines (using once more Lemma~\ref{tripot-equiv}, and Harris' characterization of the automorphisms),
no maximal tripotent generates the same class as a minmimal one.
Thus, for a spin factor $Z$, $\Delta(Z)$ has exactly $2$ elements. One of them always belongs to the spin grid, whereas
the other is in the class of the identity of $Z$ (which in the evendimensional case is equvialent to a grid element as well).
Details are in the proof of \cite[Proposition 5.1]{BoWe2}.
\end{proof}
\end{prop}

We will show that for finite dimensional $JC^*$-triple system $Z$ and $W$ we can lift $K_0$-group-homomorphisms that respect the order structure and the two scales to TRO-homomorphisms of the universal enveloping TRO of $Z$ and $W$. The additional data implemented by $\Delta(Z)$ will be used afterwards to determine the $JB^*$-isomorphisms between $Z$ and $W$.

\begin{prop}\label{lifting of TRO homs}
Let $Z$ and $W$ be finite-dimensional $JC^*$-triple systems and  $T^*(Z)=\bigoplus_{i=1}^p \mathbb{M}_{n_i,m_i}$ and
$T^*(W)=\bigoplus_{j=1}^q\mathbb{M}_{k_j,l_j}$ their universal enveloping TROs. Suppose
\begin{multline*}
\left(K_0^{\text{\tiny JB*}}(Z),K_0^{\text{\tiny JB*}}(Z)_+,\Sigma^{\text{\tiny JB*}}_\mathcal{L}(Z),\Sigma^{\text{\tiny JB*}}_\mathcal{R}(Z)\right)=\\
=\left(\Z^p,\N_0^p,\prod_{i=1}^p\{0,\ldots,n_i\},\prod_{i=1}^p\{0,\ldots,m_i\}\right)
\end{multline*}
 and
\begin{multline*}
\left(K_0^{\text{\tiny JB*}}(W),K_0^{\text{\tiny JB*}}(W)_+,\Sigma^{\text{\tiny JB*}}_\mathcal{L}(W),\Sigma^{\text{\tiny JB*}}_\mathcal{R}(W)\right)=\\
=\left(\Z^q,\N_0^q,\prod_{j=1}^q\{0,\ldots,k_j\},\prod_{j=1}^q\{0,\ldots,l_j\}\right).
\end{multline*}
Let $\alpha:K_0^{\text{\tiny JB*}}(Z)\rightarrow K_0^{\text{\tiny JB*}}(W)$ be a homomorphism such that
\begin{equation*}
\alpha\left(K_0^{\text{\tiny JB*}}(Z)_+,\Sigma^{\text{\tiny JB*}}_\mathcal{L}(Z),\Sigma^{\text{\tiny JB*}}_\mathcal{R}(Z)\right)\sub\left(K_0^{\text{\tiny JB*}}(W)_+,\Sigma^{\text{\tiny JB*}}_\mathcal{L}(W),\Sigma^{\text{\tiny JB*}}_\mathcal{R}(W)\right).
\end{equation*}
Then there exists a TRO-homomorphism $\varphi:T^*(Z)\to T^*(W)$ with
  $K_0(\varphi)=\alpha$.
\begin{proof}
We will use the fact that $\alpha:\Z^p\to\Z^q$ can be represented as a
$q\times p$-matrix $(a_{i,j})_{i,j}$ with entries $a_{i,j}\in\N_0$.
Fix $x=(z_1,\ldots z_p)\in \Sigma^{\text{\tiny JB*}}_\mathcal{L}(Z)$. Since $\alpha
(\Sigma^{\text{\tiny JB*}}_\mathcal{L}(Z))\subseteq \Sigma^{\text{\tiny JB*}}_\mathcal{L}(W)$,
\begin{equation*}
\alpha(x)=\left(\sum a_{1,i}z_i,\ldots,\sum a_{q,i}z_i\right)\leq(k_1,\ldots,k_q),
\end{equation*}
so $\sum_{j=1}^pa_{i,j}z_j\leq k_i$ for all $i=1,\ldots,q$, and similarly, $\sum_{j=1}^pa_{i,j}z_j\leq l_i$ for all $i=1,\ldots,q$.
To finish the proof, let $\varphi$ be the direct sum $\varphi:=\varphi_1\oplus\ldots\oplus\varphi_q$, where
$\varphi_j:Z\rightarrow \mathbb{M}_{k_j,l_j}$ is defined by
\[\varphi_j(x_1\oplus\ldots \oplus x_p):=\diago(\underbrace{x_1,\ldots,x_1}_{a_{1,j} times},\ldots,\underbrace{x_p,\ldots x_p}_{a_{p,j} times},0,\ldots,0),\]
for $j=1,\ldots,q$.
These TRO-homomorphisms are well-defined by the above inequalities, and $K_0(\varphi)=\alpha$.
\end{proof}
\end{prop}

\begin{thm}
Let $Z_1$ and $Z_2$ be finite-dimensional $JC^*$-triple systems. If $\sigma:K_0^{\text{ \tiny JB*}}(Z_1)\to K_0^{\text{ \tiny JB*}}(Z_2)$ is an isomorphism with $\sigma(\mathcal{KJB^*}(Z_1))=\mathcal{KJB^*}(Z_2)$, then there exists a $JB^*$-isomorphism $\vp:Z_1\to Z_2$ such that $K_0^{\text{\tiny JB*}}(\vp)=\sigma$.
\begin{proof}
We have
\begin{equation*}
\sigma\left(K_0^{\text{\tiny JB*}}(Z)_+,\Sigma^{\text{\tiny JB*}}_\mathcal{L}(Z),\Sigma^{\text{\tiny JB*}}_\mathcal{R}(Z)\right)
=\left(K_0^{\text{\tiny JB*}}(W)_+,\Sigma^{\text{\tiny JB*}}_\mathcal{L}(W),\Sigma^{\text{\tiny JB*}}_\mathcal{R}(W)\right)
\end{equation*}
and so, $\sigma$ is an isomorphism of the double-scaled ordered $K_0$-groups of $T^*(Z_1)$ and $T^*(Z_2)$. Using the above classification result for TROs, Proposition~\ref{lifting of TRO homs}, we can find a complete isometry $\vp':T^*(Z_1)\to T^*(Z_2)$ with $K_0(\vp')=\sigma$.

We will now frequently make use of the additivity of the $\mathcal{KJB^*}$-invariant, Proposition \ref{kgrid inv additiv}.
The structure of the TRO $T^*(Z_1)\simeq\bigoplus_{i=1}^p\M_{n_i,m_i}$ is obviously determined by the double-scaled ordered $K_0$-group of $T^*(Z_1)$. An inspection of the list of $\mathcal{KJB^*}$-invariants shows that
the information which is encoded in $\Delta(Z_1)$ permits to determine the Cartan-type of each summand, and
we can recover the image $\rho_{Z_1}(Z_1)\sub T^*(Z_1)$ up to (ternary) unitary equivalence (i.e an inner TRO-automorphism). The same works for $\rho_{Z_2}(Z_2)\sub T^*(Z_2)$, and it follows that $Z_1$ and $Z_2$ are isomorphic as JB*-triples.

Let $\psi:T^*(Z_2)\to T^*(Z_2)$ be the TRO isomorphism mapping $\vp'(\rho_{Z_1}(Z_1))$ to $\rho_{Z_2}(Z_2)$
(we construct $\psi$ with the aid of the universal property of $T^*(Z_2)$). Since $Z_2$ is finite dimensional
so is $T^*(Z_2)$ and thus $\psi$ is automatically unitarily equivalent to the identity, and hence the identity
on $K_0(T^*(Z_2))$. If we put
\begin{equation*}
\vp:=\rho_{Z_2}^{-1}\circ\psi\circ\vp'\circ\rho_{Z_1}:Z_1\to Z_2,
\end{equation*}
where $\rho_{Z_2}^{-1}:\rho_{Z_2}(Z_2)\to Z_2$ is the inverse of $\rho_{Z_2}$ restricted to its image,
then $\vp$ is a $JB^*$-isomorphism with $K_0^{\text{\tiny JB*}}(\vp)=\sigma$.
\end{proof}
\end{thm}

Actually, a closer look at morphisms between Cartan factors reveals the following result\cite{BoWe4}
\begin{thm}\label{JB*-homs liften}
Let $Z_1$ and $Z_2$ be finite-dimensional $JC^*$-triple systems having summands of type I, II, or III only. If $\sigma:K_0^{\text{ \tiny JB*}}(Z_1)\to K_0^{\text{ \tiny JB*}}(Z_2)$ is a homomorphism with $\sigma(\mathcal{KJB^*}(Z_1))\sub\mathcal{KJB^*}(Z_2)$, then there exists a $JB^*$-homomorphism $\vp:Z_1\to Z_2$ such that $K_0^{\text{\tiny JB*}}(\vp)=\sigma$.
\end{thm}

%

\bibliographystyle{alpha}
\bibliography{litarbeit}
\end{document}

%% file: Befehle.tex
\newcommand{\id}{\operatorname{id}}

\newcommand{\rank}{\operatorname{rank}}

\newcommand{\vp}{{\varphi}}

\newcommand{\diago}{\operatorname{diag}}

\newcommand{\Tri}{\operatorname{Tri}}

\newcommand{\R}{{\mathbb R}}

\newcommand{\cR}{{\mathcal R}}

\newcommand{\cL}{{\mathcal L}}

\newcommand{\Sp}{{\mathcal Sp}}

\newcommand{\C}{{\mathbb C}}
\newcommand{\N}{{\mathbb N}}

\newcommand{\lr}[1] {{\left(#1\right)}}
\newcommand{\fr}[1] {{\frac{1}{#1}}}

\newcommand{\set}[2]{\left\{\left.#1\right|\; #2\right\}}

\newcommand{\bt}{\begin{thm}}
\newcommand{\et}{\end{thm}}
\newcommand{\bp}{\begin{proof}}
\newcommand{\ep}{\end{proof}}
\newcommand{\bl}{\begin{lemma}}
\newcommand{\el}{\end{lemma}}
\newcommand{\bc}{\begin{cor}}
\newcommand{\ec}{\end{cor}}
\newcommand{\bd}{\begin{description}}
\newcommand{\ed}{\end{description}}

\newcommand{\fh}{\mathfrak h}

\newcommand{\ftkk}{\mathfrak{tkk}}

\newcommand{\M}{\mathbb M}
\newcommand{\Z}{\mathbb Z}
\newcommand{\Li}{\mathbb L}

\newcommand{\cH}{\mathcal H}

\newcommand{\cG}{\mathcal G}

\newcommand{\cS}{\mathcal S}

\newcommand{\sub}{\subseteq}